\documentclass[11]{article}
\usepackage[russian,english]{babel}
\usepackage{amssymb}
\usepackage{inputenc}
\usepackage{graphicx}
\begin{document}
\large
\begin{center}
\textbf{\Large On an optimal quadrature formula in Sobolev
space  $L_2^{(m)} (0,1)$}\\
Kh.M.Shadimetov, A.R.Hayotov, F.A.Nuraliev
\end{center}

\begin{abstract}
In this paper in the space $L_2^{(m)}(0,1)$ the problem of
construction of optimal quadrature formulas is considered. Here
the quadrature sum consists on values of integrand at nodes and
values of first derivative of integrand at the end points of
integration interval. The optimal coefficients are found and norm
of the error functional is calculated for arbitrary fixed $N$ and
for any $m\geq 2$. It is shown that when $m=2$ and $m=3$ the
Euler-Maclaurin quadrature formula is optimal.
\end{abstract}



\section{Introduction}

It is known, that numerical integration formulae, or quadrature
formulae, are methods for the approximate evaluation of definite
integrals. They are needed for the computation of those integrals
for which either the antiderivative of the integrand cannot be
expressed in terms of elementary functions or for which the
integrand is available only at discrete points, for example from
experimental data. In addition and even more important, quadrature
formulae provide a basic and important tool for the numerical
solution of differential and integral equations.

There are various methods in the theory of quadrature, which allow
us approximately calculate integrals with the help of finite
number values of integrand. Present work also is devoted to one of
such methods, i.e. to construction of optimal quadrature formulas
for approximate evaluation of definite integrals in the space
$L_2^{(m)}(0,1)$ equipped with the norm
$$
\|\varphi(x)\|_{L_2^{(m)}(0,1)}=\left\{\int\limits_0^1(\varphi^{(m)}(x))^2dx\right\}^{1/2}.
$$

Consider following quadrature formula
$$
\int\limits_0^1 {\varphi (x)dx \cong } \sum\limits_{\beta  = 0}^N
{C[\beta ]} \varphi [\beta ] + A\varphi '[0] + B\varphi '[N] \eqno
(1.1)
$$
with the error functional
 $$
\ell (x) = \varepsilon _{[0,1]} (x) - \sum\limits_{\beta  = 0}^N
{C[\beta ]} \delta (x - h\beta ) + A\delta '(x) + B\delta '(x - 1)
\eqno (1.2)
$$
in the space $L_2^{(m)}(0,1)$. Here $C[\beta ]$, $\beta =
\overline{0,N}$,  $A$ and $B$ are the coefficients of the formula
(1.1), $\left[\beta \right] = h\beta $, $h = \frac{1}{N},\,\,\,N =
1,2,...$ , $\varepsilon _{[0,1]} (x)$ is the indicator of interval
[0,1], $\delta (x)$ is the Dirac delta-function.

The difference
    $$
\left( {\ell (x),\varphi (x)} \right) = \int\limits_0^1 {\varphi
(x)dx - } \sum\limits_{\beta  = 0}^N {C[\beta ]} \varphi [\beta ]
- A\varphi '[0] - B\varphi '[N].
$$
is called the \emph{error} of the quadrature formula (1.1)

Error of the formula (1.1) is estimated with the help of norm of
the error functional (1.2) in the conjugate space
$L_2^{(m)*}(0,1)$, i.e. by
    $$
\left\| {\ell (x)|L_2^{(m)*} } \right\| = \mathop {\sup
}\limits_{\left\| {\varphi (x)|L_2^{(m)} } \right\| = 1} |(\ell
(x),\varphi (x))|.
$$
Furthermore, norm of the error functional  $\ell (x)$ depends on
the coefficients $C[\beta ],$  $A$ and $B$. Choice of the
coefficients when nodes are fixed is linear problem. Therefore we
minimize norm of the functional $\ell (x)$ by coefficients, i.e.
we find
    $$
\left\| {\mathop \ell \limits^ \circ  (x)|L_2^{(m)*} } \right\|  =
\mathop {\inf }\limits_{C[\beta ],A,B} \left\| {\ell
(x)|L_2^{(m)*} } \right\|. \eqno     (1.3)
$$
If $\left\| {\mathop \ell \limits^ \circ (x)|L_2^{(m)*} }
\right\|$ is found then the functional $\mathop \ell \limits^\circ
(x)$ is said to be correspond to the optimal quadrature formula
(1.1) in $L_2^{(m)} $ and corresponding coefficients are called
\emph{optimal}. Thus we get following problems.

\textbf{Problem 1.} \emph{Find norm of the error functional  $\ell
(x)$ of quadrature formula of the form (1.1) in the space
$L_2^{(m)*} (0,1)$.}

\textbf{Problem 2.} \emph{Find coefficients $C[\beta ],\,\,\,A$
and $B$ which satisfy the equality (1.3).}

Problem 2 for quadrature formulas of the form
    $$
\int\limits_0^N {\varphi (x)dx \cong } \sum\limits_{k = 0}^N {p_k
} \varphi (k)
$$
on $L_2^{(m)}$ first considered by A.Sard [1]. By A.Sard and
S.D.Meyers [2] the solution of this problem was obtained for the
following cases: $m = 1$ for arbitrary fixed $N$;  $m = 2$ for $N
\le 20$; $m = 3$ for $N \le 12$; $m = 4$ for $N \le 9$.

By I.J.Shoenberg and S.D.Silliman [3] the Sard's problem for $N
\to \infty $, i.e. for formula of the form
    $$
\int\limits_0^\infty  {\varphi (x)dx \cong \sum\limits_{k =
0}^\infty  {B_k^{(m)} \varphi (k)} }
$$
is considered. In [3] an algorithm for finding of the coefficients
$B_k^{(m)} $ is given with the help of spline of degree $2m - 1$.
In the cases $m = 2,3,...,7$ the coefficients $B_k^{(m)} $ are
calculated using a Computer. Calculation of these coefficients up
to $m = 30$ were done by  F.Ya.Zagirova [4].

In [5] in the space $L_2^{(m)} $ considered quadrature formula of
the form
$$\int\limits_{ - \eta _1 }^{N + \eta _2
} {\omega (x)\varphi (x)dx \cong \sum\limits_{\beta  = 0}^N
{C[\beta ]\varphi \left[ \beta \right]} }, \eqno  (1.4)
$$
where $0 \le \eta _j  < 1$ , $\omega (x)$ is weight function,
$C[\beta ]$ are the coefficients. In [5] the algorithm for finding
optimal coefficients $C[\beta ]$ of quadrature formulas of the
form (1.4) is given and for the optimal coefficients the system of
$2m - 2$ linear equations is obtained. These results of
S.L.Sobolev generalized above mentioned results of A.Sard,
S.D.Meyers, Schoenberg and Silliman. Further realization of
Sobolev's algorithm studied by Z.Jamalov, F.Ya.Zagirova,
Kh.M.Shadimetov. In [6], [7] the problem of construction of
optimal quadrature formulas (1.4) was completely solved for
arbitrarily fixed $N$ and for any $m$ in the space
$L_2^{(m)}(0,1)$.

Main goal of the present work is to solve problems 1 and 2 for
quadrature formulas of the form (1.1).

\section{Definitions and known formulas}

In this section we give some definitions and formulas which are
necessary in the proofs of main results.

\emph{Euler polynomials} $E_k (x)$ , $k = 1,2,...$ is defined by
following formula [8]
    $$
E_k (x) = \frac{{(1 - x)^{k + 2} }}{x}\left( {x\frac{d}{{dx}}}
\right)^k \frac{x}{{(1 - x)^2 }},\eqno     (2.1)
$$
$E_0 (x) = 1$.

For Euler polynomials following identity hold
$$
E_k (x) = x^k E_k \left( {\frac{1}{x}} \right),  \eqno (2.2)
$$
and also following theorem is take placed

\textbf{Theorem 2.1} [9]. \emph{Polynomial $P_k (x)$ which
determined by formula}
    $$
P_k (x) = (x - 1)^{k + 1} \sum\limits_{i = 0}^{k + 1}
{\frac{{\Delta ^i 0^{k + 1} }}{{(x - 1)^i }}}\eqno (2.3)
$$
\emph{is the Euler polynomial (2.1) of degree $k$, i.e. $P_k (x) =
E_k(x)$.}

Following formula is valid [10]:
    $$
\sum\limits_{\gamma  = 0}^{n - 1} {q^\gamma  \gamma ^k  =
\frac{1}{{1 - q}}\sum\limits_{i = 0}^k {\left( {\frac{q}{{1 - q}}}
\right)^i \Delta ^i 0^k  - \frac{{q^n }}{{1 - q}}\sum\limits_{i =
0}^k {\left( {\frac{q}{{1 - q}}} \right)^i \Delta ^i \gamma ^k
|_{\gamma  = n} ,} } }\eqno    (2.4)
$$
where $\Delta^i\gamma^k$ is finite difference of order $i$ of
$\gamma^k$, $\Delta^i0^k=\Delta^i\gamma^k|_{\gamma=0}$.\\
At last we give following well known formulas from [11]
    $$
\sum\limits_{\gamma  = 0}^{\beta  - 1} {\gamma ^k  =
\sum\limits_{j = 1}^{k + 1} {\frac{{k!\,B_{k + 1 - j} }}{{j!\,(k +
1 - j)!}}\,\beta ^j ,} }\eqno    (2.5)
$$
where $B_{k + 1 - j} $ are Bernoulli numbers,
    $$
\Delta ^\alpha  x^\nu   = \sum\limits_{p = 0}^\nu  {\left(
{\begin{array}{c}
   \nu   \\
   p  \\
\end{array}} \right)\Delta ^\alpha  } 0^p x^{\nu  - p}.
\eqno (2.6)
$$

\section{The extremal function and the representation of the error functional norm}

To solve problem 1, i.e. for finding norm of the error functional
(1.2) in the space $L_2^{(m)} (0,1)$ concept of the extremal
function is used [12]. The function $\psi _\ell (x)$ is said to be
\emph{extremal function} of the error functional (1.2) (see [12])
if following equality holds
    $$
\left( {\ell (x),\psi _\ell  (x)} \right) = \left\| {\ell
|L_2^{(m) * } } \right\|\left\| {\psi _\ell  |L_2^{(m)} }
\right\|.
$$
In the space $L_2^{(m)} (0,1)$ the extremal function $\psi
_\ell(x)$ of the error functional $\ell (x)$ is found by
S.L.Sobolev. This extremal function have the form
        $$
\psi _\ell  (x) = ( - 1)^m \ell (x) * G(x) + P_{m - 1} (x),\eqno
(3.1)
$$
where $G(x) = \frac{{x^{2m - 1} signx}}{{2(2m - 1)!}},$ $P_{m - 1}
(x)$ is a polynomial of degree $m - 1$. Since the functional $\ell
(x)$  belongs to the space $L_2^{(m)*}(0,1)$ therefore following
holds
    $$
(\ell (x),x^\alpha  ) = 0,\,\,\,\,\,\,\,\,\,\,\,\alpha  =
0,1,...,m - 1.\eqno    (3.2)
$$
Norm of the error functional of quadrature formula (1.1) depends
on coefficients of this formula. Indeed, since the space
$L_2^{(m)}(0,1)$ is Hilbert space, then by using  (3.1), taking
into account of Riesz theorem about common form of a linear
continuous functional on Hilbert space, we get
    $$
\left\| \ell  \right\|^2  = (\ell ,\psi _\ell  ) = ( - 1)^{m + 1}
\left[ {\frac{{A \cdot B}}{{(2m - 3)!}}} \right. - 2\left(
{A\int\limits_0^1 {\frac{{x^{2m - 2} }}{{2(2m - 2)!}}} } \right.dx
- B\left. {\int\limits_0^1 {\frac{{(x - 1)^{2m - 2} }}{{2(2m -
2)!}}} dx} \right) +
$$
$$
+ 2\sum\limits_{\beta  = 0}^N {C[\beta ]} \left( {A\frac{{(h\beta
)^{2m - 2} }}{{2(2m - 2)!}} - B\frac{{(h\beta  - 1)^{2m - 2}
}}{{2(2m - 2)!}}} \right) + 2\sum\limits_{\beta  = 0}^N {C[\beta
]} \int\limits_0^1 {\frac{{\left| {x - h\beta } \right|^{2m - 1}
}}{{2(2m - 1)!}}} dx -
$$
$$
 - \sum\limits_{\beta  = 0}^N {\sum\limits_{\gamma  = 0}^N {C[\beta ]C[\gamma ]} } \frac{{\left| {h\beta  - h\gamma } \right|^{2m - 1} }}{{2(2m - 1)!}} - \int\limits_0^1 {\int\limits_0^1 {\left. {\frac{{(x - y)^{2m - 1} sign(x - y)}}{{2(2m - 1)!}}dxdy} \right]} } .
$$
Thus, the problem 1 is solved for quadrature formulas of the form
(1.1) in the space $L_2^{(m)}(0,1)$.

\section{The system of Wiener-Hopf type}

Now we investigate problem 2. For finding of minimum of the
$\left\| \ell \right\|^2 $ under the conditions (3.2) Lagrange
method of undetermined multipliers is used. For this we consider
following function
    $$
\Psi  = \left\| \ell  \right\|^2  + 2 \cdot ( - 1)^{m + 1}
\sum\limits_{\alpha  = 0}^{m - 1} {\lambda _\alpha  } (\ell
,x^\alpha  ).
$$
Equating to zero partial derivatives by coefficients $C[\beta
],\,\,A$ and $B,$  together with conditions (3.2) we get following
system of linear equations
                $$
\sum\limits_{\gamma  = 0}^N {C[\gamma ]} \frac{{\left| {h\beta  -
h\gamma } \right|^{2m - 1} }}{{2(2m - 1)!}} - A\frac{{(h\beta
)^{2m - 2} }}{{2(2m - 2)!}} +
$$
 $$
+ B\frac{{(h\beta  - 1)^{2m - 2} }}{{2(2m - 2)!}} +
\sum\limits_{\alpha  = 0}^{m - 1} {\lambda _\alpha  (h\beta
)^\alpha   = \int\limits_0^1 {\frac{{\left| {x - h\beta }
\right|^{2m - 1} }}{{2(2m - 1)!}}} } dx,\,\,\,\,\,\,\beta  =
\overline {0,N} , \eqno   (4.1)
$$
$$
\sum\limits_{\gamma  = 0}^N {C[\gamma ]} \frac{{(h\gamma )^{2m -
2} }}{{2(2m - 2)!}} + \frac{B}{{2(2m - 3)!}} - \lambda _1  =
\frac{1}{{2(2m - 1)!}},\eqno     (4.2)
$$
$$
\sum\limits_{\gamma  = 0}^N {C[\gamma ]} \frac{{(h\gamma  - 1)^{2m
- 2} }}{{2(2m - 2)!}} - \frac{A}{{2(2m - 3)!}} +
\sum\limits_{\alpha  = 1}^{m - 1} {\alpha \lambda _\alpha  }  =
\frac{1}{{2(2m - 1)!}},\eqno    (4.3)
$$
        $$
\sum\limits_{\gamma  = 0}^N {C[\gamma ]}  = 1,\eqno (4.4)
$$
                 $$
\sum\limits_{\gamma  = 0}^N {C[\gamma ]h\gamma  + A + B =
\frac{1}{2},}\eqno  (4.5)
$$
                 $$
\sum\limits_{\gamma  = 0}^N {C[\gamma ](h\gamma )^\alpha   +
\alpha B = \frac{1}{{\alpha  + 1}},} \,\,\,\,\,\,\alpha  =
\overline {2,m - 1}.\eqno     (4.6)
$$
The system (4.1)-(4.6) is called by \emph{system of Wiener-Hof
type} for the optimal coefficients [12]. In the system (4.1)-(4.6)
coefficients $C[\beta ],$ $\beta = \overline {0,N}$, $A$ and $B,$
and also $\lambda _\alpha ,\,\,\,\,\alpha  = \overline {0,m - 1}$
are unknowns. The system (4.1)-(4.6) has unique solution. The
proof of existence and uniqueness of the solution of this system
is as the proof of existence and uniqueness of the solution of
Wiener-Hopf system of the optimal coefficients in the space
$L_2^{(m)}(0,1)$ for quadrature formulas of the form (1.4) (see
[13]).

\section{The optimal coefficients and norm of the error functional}

In present section we study solution of the system (4.1)-(4.6). In
the solution of this system we use the approach which used in
solution of the linear system for optimal coefficients of
quadrature formulas of the form (1.4) in [6].

\subsection{The coefficients of optimal quadrature formulas}

It is easy to prove following lemma for the coefficients  $C[\beta
]$ of quadrature formulas of the form (1.1).

\textbf{Lemma 5.1.} \textit{The optimal coefficients $C[\beta ]$,
$1 \le \beta \le N - 1$, of quadrature formulas of the form (1.1)
have following form
    $$
C[\beta ] = h\left( {1 + \sum\limits_{k = 1}^{m - 1} {\left( {d_k
q_k^\beta   + p_k q_k^{N - \beta } } \right)} }
\right),\,\,\,\,\,\,\,\,\,\,\,\,\,\,1 \le \beta  \le N - 1,\eqno
(5.1)
$$
where $d_{k,} \,p_k  $ are unknowns, $q_k  $ are roots of the
Euler polynomial $E_{2m - 2} (q),\,\,\,\left| {q_k } \right| < 1.$
}

Lemma is proved as lemma 3 of the work [9] and in the proof the
discrete analogue $D_m [\beta ]$ of the polyharmonic operator
$\frac{d^{2m}}{dx^{2m}}$ is used. The discrete analogue $D_m
[\beta ]$ of the polyharmonic operator $\frac{d^{2m}}{dx^{2m}}$ is
constructed in [14].

We need following lemmas in proof of main results.

\textbf{Lemma 5.2.} \textit{Following identity is take placed
    $$
\sum\limits_{i = 0}^\alpha  {\frac{{dq + pq_{}^{N + i} ( - 1)^{i +
1} }}{{(q - 1)^{i + 1} }}\Delta ^i 0^\alpha  }  = ( - 1)^{\alpha +
1} \sum\limits_{i = 0}^\alpha  {\frac{{dq_{}^i  + pq_{}^{N + 1} (
- 1)^{i + 1} }}{{(1 - q)^{i + 1} }}\Delta ^i 0^\alpha  },\eqno
(5.2)
$$
here $\alpha $ and $N$ are natural numbers, $\Delta ^i 0^\alpha$
is finite difference of order $i$ of $\gamma ^\alpha $ at the
point 0.}

\textbf{Proof.}  For convenience left and right sides of (5.2)we
denote by $L_1 $ и $L_2$ respectively, i.e.
    $$
L_1  = \sum\limits_{i = 0}^\alpha  {\frac{{dq + pq_{}^{N + i} ( -
1)^{i + 1} }}{{(q_{}  - 1)^{i + 1} }}\Delta ^i 0^\alpha  } \mbox{
и }L_2  = ( - 1)^{\alpha  + 1} \sum\limits_{i = 0}^\alpha
{\frac{{dq_{}^i  + pq_{}^{N + 1} ( - 1)^{i + 1} }}{{(1 - q_{} )^{i
+ 1} }}\Delta ^i 0^\alpha  } .$$

First consider $L_1 $. By using the equality (2.3) and identity
 (2.2) for $L_1 $ consequently we get
$$
L_1  = \sum\limits_{i = 0}^\alpha  {\frac{{dq + pq_{}^{N + i} ( -
1)^{i + 1} }}{{(q_{}  - 1)^{i + 1} }}\Delta ^i 0^\alpha  }  =
\frac{{dq}}{{(q - 1)^{\alpha  + 1} }}E_{\alpha  - 1} (q) +
\frac{{pq_{}^{N + \alpha } ( - 1)^{\alpha  + 1} }}{{(q -
1)^{\alpha  + 1} }}E_{\alpha  - 1} \left( {\frac{1}{q}} \right) =
$$
$$
= \frac{{dq}}{{(q - 1)^{\alpha  + 1} }}E_{\alpha  - 1} (q) +
\frac{{pq_{}^{N + \alpha } ( - 1)^{\alpha  + 1} }}{{(q -
1)^{\alpha  + 1} }}\frac{{E_{\alpha  - 1} (q)}}{{q_{}^{\alpha  -
1} }} = \frac{{dq + pq_{}^{N + 1} ( - 1)^{\alpha  + 1} }}{{(q -
1)^{\alpha  + 1} }}E_{\alpha  - 1} (q).\eqno (5.3)
$$
Similarly for $L_2 $ by using  (2.3) and (2.2) we have
$$
L_2  = \sum\limits_{i = 0}^\alpha  {\frac{{dq_{}^i  + pq_{}^{N +
1} ( - 1)^{i + 1} }}{{(1 - q)^{i + 1} }}\Delta ^i 0^\alpha  }  =
\frac{{dq_{}^\alpha  }}{{(q - 1)^{\alpha  + 1} }}E_{\alpha  - 1}
\left( {\frac{1}{q}} \right) + \frac{{pq_{}^{N + 1} }}{{(q -
1)^{\alpha  + 1} }}E_{\alpha  - 1} \left( q \right) =
$$
$$
= \frac{{dq}}{{(1 - q)^{\alpha  + 1} }}E_{\alpha  - 1} (q) +
\frac{{pq_{}^{N + 1} }}{{(q - 1)^{\alpha  + 1} }}E_{\alpha  - 1}
(q) = \frac{{dq( - 1)^{\alpha  + 1}  + pq_{}^{N + 1} }}{{(q -
1)^{\alpha  + 1} }}E_{\alpha  - 1} (q) =
$$
$$
= ( - 1)^{\alpha  + 1} \frac{{dq + pq_{}^{N + 1} ( - 1)^{\alpha  +
1} }}{{(q - 1)^{\alpha  + 1} }}E_{\alpha  - 1} (q).\eqno
(5.4)
$$
From (5.3) and (5.4) clear, that $L_1  = ( - 1)^{\alpha  + 1} L_2
$. Lemma 5.2 is proved.

We denote
    $$
Z_p  = \sum\limits_{k = 1}^{m - 1} {\sum\limits_{i = 0}^p
{\frac{{d_k q_k^{N + i}  + p_k q_k ( - 1)^{i + 1} }}{{(1 - q_k
)^{i + 1} }}\Delta ^i 0^p } }. \eqno  (5.4*)
$$

\textbf{Lemma 5.3.} \emph{Following identities are valid}
$$
\sum\limits_{j = 1}^{m - 1} {\frac{{( - 1)^{j - 1} }}{{(j -
1)!}}\sum\limits_{i = 1}^{2m - 3 - j} {\frac{{B_{2m - j - i} h^{2m
- j - i} }}{{i!\,\,(2m - j - i)!}} = } }
$$
 $$
= \sum\limits_{j = 3}^m {\frac{{B_j h^j }}{{j!}}\sum\limits_{i =
0}^{m - 2} {\frac{{( - 1)^i }}{{i!\,\,(2m - 1 - j - i)!}} + } }
\sum\limits_{j = m + 1}^{2m - 2} {\frac{{B_j h^j
}}{{j!}}\sum\limits_{i = 0}^{2m - 2 - j} {\frac{{( - 1)^i
}}{{i!\,\,(2m - 1 - j - i)!}}} }
$$
and
$$
\sum\limits_{j = 1}^{m - 1} {\frac{{( - 1)^{j - 1} }}{{(j -
1)!}}\sum\limits_{p = 2}^{2m - 1 - j} {\frac{{h^{p + 1} Z_p
}}{{p!\,\,(2m - 1 - j - p)!}} = } }
$$
$$
= \sum\limits_{j = 3}^{m + 1} {\frac{{h^j Z_{j - 1} }}{{(j -
1)!}}\sum\limits_{l = 0}^{m - 2} {\frac{{( - 1)^l }}{{l!\,\,(2m -
1 - j - l)!}} + } } \sum\limits_{j = m + 2}^{2m - 1} {\frac{{h^j
Z_{j - 1} }}{{(j - 1)!}}\sum\limits_{l = 0}^{2m - 1 - j} {\frac{{(
- 1)^l }}{{l!\,\,(2m - 1 - j - l)!}}} }.
$$

The proof of lemma is obtained by expansion in powers of $h$ of
left sides of given identities.

For the coefficients of optimal quadrature formulas of the form
(1.1) following theorem holds.

\textbf{Theorem 5.1.} \textit{Among quadrature formulas of the
form  (1.1) with the error functional (1.2) there exists unique
optimal formula which coefficients are determined by following
formulas
    $$
C[0] = h\left( {\frac{1}{2} + \sum\limits_{k = 1}^{m - 1}
{\frac{{p_k q_k^N  - d_k q_k }}{{1 - q_k }}} } \right),\eqno
(5.5)
$$
  $$
C[\beta ] = h\left( {1 + \sum\limits_{k = 1}^{m - 1} {(d_k
q_k^\beta   + p_k q_k^{N - \beta } )} }
\right),\,\,\,\,\,\,\,\,\,\,\,\,\,\,\,\beta  = \overline {1,N -
1},\eqno   (5.6)
$$
$$
C[N] = h\left( {\frac{1}{2} + \sum\limits_{k = 1}^{m - 1}
{\frac{{d_k q_k^N  - p_k q_k }}{{1 - q_k }}} } \right),\eqno
(5.7)
$$
$$
A = h^2 \left( {\frac{1}{{12}} - \sum\limits_{k = 1}^{m - 1}
{\frac{{d_k q_k  + p_k q_k^{N + 1} }}{{(1 - q_k )^2 }}} }
\right),\eqno     (5.8)
$$
$$
B = h^2 \left( { - \frac{1}{{12}} + \sum\limits_{k = 1}^{m - 1}
{\frac{{d_k q_k^{N + 1}  + p_k q_k }}{{(1 - q_k )^2 }}} }
\right),\eqno     (5.9)
$$
where $d_k $ and $p_k $ satisfy following system $2m - 2$ linear
equations:
                $$
\sum\limits_{k = 1}^{m - 1} {\sum\limits_{i = 0}^j {\frac{{d_k q_k
+ p_k q_k^{N + i} ( - 1)^{i + 1} }}{{(q_k  - 1)^{i + 1} }}} }
\Delta ^i 0^j  = \frac{{B_{j + 1} }}{{j + 1}},\,\,\,\,\,\,j =
\overline {2,m - 1},\eqno     (5.10)
$$
$$
\sum\limits_{k = 1}^{m - 1} {\sum\limits_{i = 0}^{2m - 2}
{\frac{{d_k q_k  + p_k q_k^{N + i} ( - 1)^{i + 1} }}{{(q_k  -
1)^{i + 1} }}} } \Delta ^i 0^{2m - 2}  = 0,\eqno       (5.11)
$$
$$
\sum\limits_{k = 1}^{m - 1} {\sum\limits_{i = 0}^j {(1 - q_k^N
)\frac{{( - 1)^{i + 1} d_k q_k^i  - p_k q_k }}{{(q_k  - 1)^{i + 1}
}}} } \Delta ^i 0^j  = 0,\,\,\,\,\,\,\,j = \,\overline {2,m -
1}.\eqno     (5.12)
$$
$$
\sum\limits_{k = 1}^{m - 1} {\sum\limits_{i = 0}^{2m - 2}
{\frac{{( - 1)^{i + 1} d_k q_k^{N + i}  + p_k q_k^{} }}{{(q_k  -
1)^{i + 1} }}} } \Delta ^i 0^{2m - 2}  = 0.\eqno      (5.13)
$$
Here $B_\alpha $ are Bernoulli numbers, $\Delta ^i \gamma ^j  $ is
difference of order $i$ of $\gamma^j $, $\,\Delta ^i 0^j = \Delta
^i \gamma ^j |_{\gamma  = 0} $ , $\,q_k $ are roots of Euler
polynomial of degree $2m - 2$ ,$\,\left| {q_k } \right| < 1$. }

\textbf{Proof.} First we give plan of proof.

From (5.1) clear that instead of unknowns $C[\beta ]$, $\beta =
\overline {1,N - 1} $ it is sufficient to find unknowns $d_k $,
$p_k $, $k = \overline {1,m - 1} $. The coefficients $C[0]$,
$C[N]$, $A$, $B$ and $\lambda _\alpha ,\,\,\,\alpha  = \overline
{0,m - 1}$ are expressed by $d_k$ and $p_k $, $k = \overline {1,m
- 1}$. So if we find $d_k $ and $p_k $, then the system
(4.1)-(4.6) is solved completely. Substituting the equality (5.1)
to equation (4.1) we get polynomial of degree $2m$ of $h\beta $ on
both sides of (4.1). Equating coefficients of same degrees of
$h\beta $ we find $\lambda _\alpha ,\,\,\,\alpha  = \overline {0,m
- 1} $, $C[0]$, $A$ and system (5.10) for $d_k $, $p_k $. Taking
account of (5.1), (5.5), (5.9), from conditions (4.4) and (4.5) we
get (5.7), (5.9), i.e. we obtain $C[N]$ and $B$. Further, by using
(5.1), (5.9) and expression for $\lambda _1 $, from (4.2) we get
the equation (5.11). System of equations (5.12) for unknowns $d_k
$ ,$p_k $ , we obtain from equation (4.6), using (5.1),
(5.5)-(5.9). Finally, taking into account (5.1), (5.8) and
$\lambda _\alpha ,\,\,\,\alpha = \overline {1,m - 1}$ from
equation (4.3) we have (5.13).

Further we give detailed explanation of proof of the theorem.

First we consider  first sum of equation (4.1). For this sum we
have
$$
S = \sum\limits_{\gamma  = 0}^N {C[\gamma ]\frac{{|h\beta  -
h\gamma |^{2m - 1} }}{{2(2m - 1)!}} = }
$$
$$=C[0]\frac{{(h\beta )^{2m - 1} }}{{(2m - 1)!}} +
\sum\limits_{\gamma  = 1}^\beta  {C[\gamma ]\frac{{(h\beta  -
h\gamma )^{2m - 1} }}{{(2m - 1)!}}}  - \sum\limits_{\gamma  = 0}^N
{C[\gamma ]} \frac{{(h\beta  - h\gamma )^{2m - 1} }}{{2(2m -
1)!}}.
$$
Lat two sums of the expression $S$ we denote
$$
S_1  = \sum\limits_{\gamma  = 1}^\beta  {C[\gamma ]\frac{{(h\beta
- h\gamma )^{2m - 1} }}{{(2m - 1)!}}},\ \ \  S_2  =
\sum\limits_{\gamma  = 0}^N {C[\gamma ]} \frac{{(h\beta  - h\gamma
)^{2m - 1} }}{{2(2m - 1)!}}
$$
and we calculate them separately.\\
By using lemma 5.1 and formulas (2.4), (2.5) for $S_1 $ we have
$$
S_1  = \sum\limits_{\gamma  = 0}^\beta
{h\left( {1 + \sum\limits_{k = 1}^{m - 1} {\left( {d_k q_k^\gamma
+ p_k q_k^{N - \gamma } } \right)} } \right)\frac{{(h\beta  -
h\gamma )^{2m - 1} }}{{(2m - 1)!}} = } $$
$$ = \frac{{h^{2m}
}}{{(2m - 1)!}}\left[ {\sum\limits_{\gamma  = 0}^{\beta  - 1}
{\gamma ^{2m - 1}  + \sum\limits_{k = 1}^{m - 1} {\left( {d_k
q_k^\beta \sum\limits_{\gamma  = 0}^{\beta  - 1} {q_k^{ - \gamma }
\gamma ^{2m - 1} }  + p_k q_k^{N - \beta } \sum\limits_{\gamma  =
0}^{\beta  - 1} {q_k^\gamma  \gamma ^{2m - 1} } } \right)} } }
\right] = $$ $$ = \frac{{h^{2m} }}{{(2m - 1)!}}\left[
{\sum\limits_{j = 1}^{2m} {\frac{{(2m - 1)!B_{2m - j} }}{{j! \cdot
(2m - j)!}}\beta ^j }  + } \right.\sum\limits_{k = 1}^{m - 1}
{\left[ {d_k q_k^\beta \left\{ {\frac{{q_k }}{{q_k  -
1}}\sum\limits_{i = 0}^{2m - 1} {\frac{{\Delta ^i 0^{2m - 1}
}}{{(q_k  - 1)^i }}} } \right.} \right.}  -
$$
$$
- \frac{{q_k^{1 - \beta } }}{{q_k  - 1}}\left. {\sum\limits_{i =
0}^{2m - 1} {\frac{{\Delta ^i \beta ^{2m - 1} }}{{(q_k  - 1)^i }}}
} \right\} + p_k q_k^{N - \beta } \left\{ {\frac{1}{{1 - q_k
}}\sum\limits_{i = 0}^{2m - 1} {\left( {\frac{{q_k }}{{q_k  - 1}}}
\right)^i \Delta ^i 0^{2m - 1}  - } } \right.$$ $$ -
\frac{{q_k^\beta  }}{{1 - q_k }}\left. {\left. {\left.
{\sum\limits_{i = 0}^{2m - 1} {\left( {\frac{{q_k }}{{q_k  - 1}}}
\right)^i \Delta ^i \beta ^{2m - 1} } } \right\}} \right]}
\right].$$ Taking into account that $q_k $ is a root of Eulaer
polynomial $E_{2m - 2} (q)$ an using formulas (2.3), (2.6) the
expression for $S_1$ we reduce to following form
$$S_1  = \frac{{(h\beta )^{2m} }}{{(2m)!}} +
h \cdot \frac{{(h\beta )^{2m - 1} }}{{(2m - 1)!}}B_1  + h^{2m}
\sum\limits_{j = 1}^{2m - 2} {\frac{{B_{2m - j} }}{{j!(2m -
j)!}}\beta ^j  + } $$ $$ + h^{2m} \sum\limits_{j = 0}^{2m - 1}
{\frac{{\beta ^{2m - 1 - j} }}{{j!(2m - 1 - j)!}}\sum\limits_{k =
1}^{m - 1} {\sum\limits_{i = 0}^j {\frac{{ - d_k q_k  + p_k q_k^{N
+ i} ( - 1)^i }}{{(q_k  - 1)^{i + 1} }}\Delta ^i 0^j .} } }
\eqno            (5.14)$$ Now consider $S_2 $. By using conditions
of orthogonality (4.4)-(4.6) the expression $S_2$ we rewrite by
powers of $h\beta$
$$S_2  = \sum\limits_{\gamma  = 0}^N C[\gamma ]\frac{{(h\beta  -
h\gamma )^{2m - 1} }}{{2(2m - 1)!}} = \frac{1}{2}\sum\limits_{j =
2}^{m - 1} \frac{{(h\beta )^{2m - 1 - j} }}{{j!(2m - 1 -
j)!}}\left( {\frac{1}{{j + 1}} - jB} \right) -
$$
$$
-\frac{{(h\beta )^{2m - 2} }}{{2(2m - 2)!}}\left( {\frac{1}{2} - A
- B} \right) + \frac{{(h\beta )^{2m - 1} }}{{2(2m - 1)!}} +
\frac{1}{2}\sum\limits_{j = m}^{2m - 1} {\frac{{(h\beta )^{2m - 1
- j} }}{{j!(2m - 1 - j)!}}\sum\limits_{\gamma  = 0}^N {C[\gamma ](
- h\gamma )^j .} } \eqno (5.15)
$$
The right side of the equation (4.1) have following form
    $$
\int\limits_0^1 {\frac{{|x - h\beta |^{2m - 1} }}{{2(2m - 1)!}}dx
= \frac{{(h\beta )^{2m} }}{{(2m)!}} + \sum\limits_{j = 0}^{2m - 1}
{\frac{{( - h\beta )^{2m - 1 - j} }}{{2(2m - 1 - j)!(j + 1)!}}} }
\eqno    (5.16)$$ Substituting (5.16) and $S$ into equation (4.1)
and using (5.14), (5.15) we have
$$\frac{{(h\beta )^{2m} }}{{(2m)!}} + C[0]\frac{{(h\beta )^{2m - 1}
}}{{(2m - 1)!}} + h\frac{{(h\beta )^{2m - 1} }}{{(2m - 1)!}}B_1  +
\sum\limits_{j = 1}^{2m - 2} {\frac{{B_{2m - j} h^{2m - j} (h\beta
)^j }}{{j!(2m - j)!}} + } $$
$$
+ \sum\limits_{j = 0}^{2m - 1} {\frac{{h^{j + 1} (h\beta )^{2m - 1
- j} }}{{j!(2m - 1 - j)!}}\sum\limits_{k = 1}^{m - 1}
{\sum\limits_{i = 0}^j {\frac{{ - d_k q_k  + p_k q_k^{N + i} ( -
1)^i }}{{(q_k  - 1)^{i + 1} }}\Delta ^i 0^j  - } } }
$$
$$
- \frac{1}{2}\sum\limits_{j = 2}^{m - 1} {\frac{{(h\beta )^{2m - 1
- j} ( - 1)^j }}{{j!(2m - 1 - j)!}}\left( {\frac{1}{{j + 1}} - jB}
\right) + \frac{{(h\beta )^{2m - 2} }}{{2\,\,(2m - 2)!}}\left(
{\frac{1}{2} - A - B} \right) - }
$$
$$
- \frac{{(h\beta )^{2m - 1} }}{{2(2m - 1)!}} -
\frac{1}{2}\sum\limits_{j = m}^{2m - 1} {\frac{{(h\beta )^{2m - 1
- j} }}{{j!(2m - 1 - j)!}}} \sum\limits_{\gamma  = 0}^N {C[\gamma
]} ( - h\gamma )^j  - A\frac{{(h\beta )^{2m - 2} }}{{2(2m - 2)!}}
+
$$
$$
+ B\sum\limits_{j = 0}^{2m - 2} {\frac{{(h\beta )^{2m - 2 - j} ( -
1)^j }}{{2 \cdot j! \cdot (2m - 2 - j)!}}}  + \sum\limits_{\alpha
= 0}^{m - 1} {\lambda _\alpha  (h\beta )^\alpha  }  =
\frac{{(h\beta )^{2m} }}{{(2m)!}} + \sum\limits_{j = 0}^{2m - 1}
{\frac{{( - h\beta )^{2m - 1 - j} }}{{2 \cdot (2m - 1 - j)! \cdot
(j + 1)!}}}.
$$
Hence equating coefficients of same powers of $h\beta $ gives
$$\lambda _j  = \frac{1}{{(2m - 1 - j)! \cdot j!}}\left( {\frac{{(
- 1)^{2m - j} }}{{2(2m - j)}} - \frac{{B_{2m - j} h^{2m - j}
}}{{2m - j}}} \right.-
$$
$$
- h^{2m - j} \sum\limits_{k = 1}^{m - 1} {\sum\limits_{i = 0}^{2m
- 1 - j} {\frac{{ - d_k q_k  + p_k q_k^{N + i} ( - 1)^i }}{{(q_k
- 1)^{i + 1} }}\Delta ^i 0^{2m - 1 - j}  + } } $$
$$ + \frac{1}{2}\sum\limits_{\gamma  = 0}^N {C[\gamma ]} ( -
h\gamma )^{2m - 1 - j}  - (2m - 1 - j) \cdot B \cdot \left.
{\frac{{( - 1)^{2m - 2 - j} }}{2}} \right),\,\,\,\,\,j = 1,2,...,m
- 1,\eqno (5.17)$$
$$
\lambda _0  = \frac{1}{{2 \cdot (2m)!}} + \frac{1}{{2 \cdot (2m -
1)!}}\sum\limits_{\gamma  = 0}^N {C[\gamma ]} ( - h\gamma )^{2m -
1}  - B\frac{1}{{2 \cdot (2m - 2)!}}, \eqno  (5.18)
$$
$$
\sum\limits_{k = 1}^{m - 1} {\sum\limits_{i = 0}^j {\frac{{d_k q_k
+ p_k q_k^{N + i} ( - 1)^{i + 1} }}{{(q_k  - 1)^{i + 1} }}} }
\Delta ^i 0^j  = \frac{{B_{j + 1} }}{{j + 1}},\,\,\,\,\,j =
\overline {2,m - 1} ,\eqno         (5.19)
$$
 $$
C[0] = h\left( {\frac{1}{2} + \sum\limits_{k = 1}^{m - 1}
{\frac{{p_k q_k^N  - d_k q_k }}{{1 - q_k }}} } \right),\eqno
(5.20)
$$
 $$
A = h^2 \left( {\frac{1}{{12}} - \sum\limits_{k = 1}^{m - 1}
{\frac{{d_k q_k^{}  + p_k q_k^{N + 1} }}{{(1 - q_k )^2 }}} }
\right). \eqno  (5.21)
$$
Here equation (5.19) is the equation (5.10) for unknowns $d_k $
and $p_k $.

Substituting expressions (5.20) and (5.21)  into (4.4) and (4.5),
also taking into account (5.1), we find $C[N]$ and $B$, which have
following form
    $$
C[N] = h\left( {\frac{1}{2} + \sum\limits_{k = 1}^{m - 1}
{\frac{{d_k q_k^N  - p_k q_k }}{{1 - q_k }}} } \right),\eqno
(5.22)
$$
$$
B = h^2 \left( { - \frac{1}{{12}} + \sum\limits_{k = 1}^{m - 1}
{\frac{{d_k q_k^{N + 1}+p_k q_k }}{{(1 - q_k )^2 }}} } \right)
.\eqno   (5.23)
$$
Now substituting the expression of $\lambda _1 $  from (5.17) when
$j = 1$ into  (4.2), we get one more equation with respect to
unknowns  $d_k $ and $p_k $:
    $$
\sum\limits_{k = 1}^{m - 1} {\sum\limits_{i = 0}^{2m - 2}
{\frac{{d_k q_k  + p_k q_k^{N + i} ( - 1)^{i + 1} }}{{(q_k  -
1)^{i + 1} }}} } \Delta ^i 0^{2m - 2}  = 0,\eqno   (5.24)
$$
i.e. we obtain the equation (5.11).

Next, to obtain (5.12) we use equation (4.6). Consider equation
(4.6). Since $\alpha  = \overline {2,m - 1} $, then
    $$
\sum\limits_{\gamma  = 0}^N {C[\gamma ](h\gamma )^\alpha   +
\alpha \,B = } \sum\limits_{\gamma  = 1}^{N - 1} {C[\gamma
](h\gamma )^\alpha   + C[N] + \alpha \,B = } \frac{1}{{\alpha  +
1}}.\eqno     (5.25)
$$
We denote $L = \sum\limits_{\gamma  = 1}^{N - 1} {C[\gamma
](h\gamma )^\alpha  } $. Using (5.1), (2.4), (2.5) for $L$ we get
    $$
L = h^{\alpha  + 1} \left( {\sum\limits_{\gamma  = 1}^{N - 1}
{\gamma ^\alpha  }  + \sum\limits_{k = 1}^{m - 1} {\left( {d_k
\sum\limits_{\gamma  = 1}^{N - 1} {q_k^\gamma  \gamma ^\alpha  }
+ p_k q_k^N \sum\limits_{\gamma  = 1}^{N - 1} {q_k^{ - \gamma }
\gamma ^\alpha  } } \right)} } \right) =
$$
$$
= \sum\limits_{j = 1}^{\alpha  + 1} {\frac{{\alpha !B_{\alpha  + 1
- j} }}{{j!\,(\alpha  + 1 - j)!}}} h^{\alpha  + 1 - j}  +
h^{\alpha  + 1} \sum\limits_{k = 1}^{m - 1} {\sum\limits_{i =
0}^\alpha  {\frac{{d_k q_k^i  + p_k q_k^{N + 1} ( - 1)^{i + 1}
}}{{(1 - q_k )^{i + 1} }}} \Delta ^i 0^\alpha  }  -
$$
$$
- h^{\alpha  + 1} \sum\limits_{k = 1}^{m - 1} {\sum\limits_{i =
0}^\alpha  {\frac{{d_k q_k^{N + i}  + p_k q_k^{} ( - 1)^{i + 1}
}}{{(1 - q_k )^{i + 1} }}} \Delta ^i N^\alpha  }.
$$
Hence taking into account(2.6) and grouping in powers of $h$, we
have
    $$
L = \frac{1}{{\alpha  + 1}} + \sum\limits_{j = 1}^\alpha
{\frac{{\alpha !h^j }}{{(j - 1)!\,\,(\alpha  + 1 - j)!}}} \left(
{\frac{{B_j }}{j} - \sum\limits_{k = 1}^{m - 1} {\sum\limits_{i =
0}^{j - 1} {\frac{{d_k q_k^{N + i}  + p_k q_k ( - 1)^{i + 1}
}}{{(1 - q_k )^{i + 1} }}} \Delta ^i 0^{j - 1} } } \right) +
$$
$$
+ h^{\alpha  + 1} \left( {\sum\limits_{k = 1}^{m - 1}
{\sum\limits_{i = 0}^\alpha  {\frac{{d_k q_k^i  + p_k q_k^{N + 1}
( - 1)^{i + 1} }}{{(1 - q_k )^{i + 1} }}} \Delta ^i 0^\alpha   -
\sum\limits_{k = 1}^{m - 1} {\sum\limits_{i = 0}^\alpha
{\frac{{d_k q_k^{N + i}  + p_k q_k ( - 1)^{i + 1} }}{{(1 - q_k
)^{i + 1} }}} \Delta ^i 0^\alpha  } } } \right).
$$
Substitution of obtained expression of $L$ to (5.25) gives
 $$
\frac{1}{{\alpha  + 1}} + \sum\limits_{j = 1}^\alpha
{\frac{{\alpha !h^j }}{{(j - 1)!\,\,(\alpha  + 1 - j)!}}} \left(
{\frac{{B_j }}{j} - \sum\limits_{k = 1}^{m - 1} {\sum\limits_{i =
0}^{j - 1} {\frac{{d_k q_k^{N + i}  + p_k q_k ( - 1)^{i + 1}
}}{{(1 - q_k )^{i + 1} }}} \Delta ^i 0^{j - 1} } } \right) +
$$
$$
+ h^{\alpha  + 1} \left( {\sum\limits_{k = 1}^{m - 1}
{\sum\limits_{i = 0}^\alpha  {(1 - q_k^N )\frac{{( - 1)^{i + 1}
d_k q_k^i  - p_k q_k^{} }}{{(q_k  - 1)^{i + 1} }}} \Delta ^i
0^\alpha  } } \right) + C[N] + \alpha B = \frac{1}{{\alpha  + 1}}.
$$
Hence keeping in mind (5.22), (5.23) we get
    $$
\sum\limits_{j = 3}^\alpha  {\frac{{\alpha !h^j }}{{(j -
1)!\,\,(\alpha  + 1 - j)!}}} \left( {\frac{{B_j }}{j} -
\sum\limits_{k = 1}^{m - 1} {\sum\limits_{i = 0}^{j - 1}
{\frac{{d_k q_k^{N + i}  + p_k q_k ( - 1)^{i + 1} }}{{(1 - q_k
)^{i + 1} }}} \Delta ^i 0^{j - 1} } } \right) +
$$
$$
+ h^{\alpha  + 1} \left( {\sum\limits_{k = 1}^{m - 1}
{\sum\limits_{i = 0}^\alpha  {(1 - q_k^N )\frac{{( - 1)^{i + 1}
d_k q_k^i  - p_k q_k^{} }}{{(q_k  - 1)^{i + 1} }}} \Delta ^i
0^\alpha  } } \right) = 0.\eqno             (5.26)
$$
Clearly that the left side of (5.26) is polynomial of degree
$\alpha+ 1$ with respect to $h$. From (5.26) we obtain that each
coefficient of this polynomial is equal to zero, i.e.
    $$
\sum\limits_{k = 1}^{m - 1} {\sum\limits_{i = 0}^{j - 1}
{\frac{{d_k q_k^{N + i}  + p_k q_k ( - 1)^{i + 1} }}{{(1 - q_k
)^{i + 1} }}} \Delta ^i 0^{j - 1} }  = \frac{{B_j
}}{j},\,\,\,\,\,j = \overline {3,\alpha }  \eqno  (5.27)
$$
and
$$
\sum\limits_{k = 1}^{m - 1} {\sum\limits_{i = 0}^\alpha  {(1 -
q_k^N )\frac{{( - 1)^{i + 1} d_k q_k^i  - p_k q_k^{} }}{{(q_k  -
1)^{i + 1} }}} \Delta ^i 0^\alpha  }  = 0.          \eqno (5.28)
$$
Since  $\alpha  = \overline {2,m - 1} $ then from (5.27) and
(5.28) we have
    $$
\sum\limits_{k = 1}^{m - 1} {\sum\limits_{i = 0}^j {\frac{{d_k
q_k^{N + i}  + p_k q_k ( - 1)^{i + 1} }}{{(1 - q_k )^{i + 1} }}}
\Delta ^i 0^j }  = \frac{{B_{j + 1} }}{{j + 1}},\,\,\,\,\,j =
\overline {2,m - 2}\eqno    (5.29)
$$
and
$$
\sum\limits_{k = 1}^{m - 1} {\sum\limits_{i = 0}^\alpha
{\frac{{d_k q_k^i  + p_k q_k^{N + 1} ( - 1)^{i + 1} }}{{(1 - q_k
)^{i + 1} }}} \Delta ^i 0^\alpha   = \sum\limits_{k = 1}^{m - 1}
{\sum\limits_{i = 0}^\alpha  {\frac{{d_k q_k^{N + i}  + p_k q_k (
- 1)^{i + 1} }}{{(1 - q_k )^{i + 1} }}} \Delta ^i 0^\alpha
,\,\,\,\,\alpha  = \overline {2,m - 1} } }. \eqno  (5.30)
$$
Hence by using (5.30) and lemma 5.2 it is easy to show the system
of equations (5.29) is the part of (5.19). Thus here we get only
system of equations (5.30) which is (5.12).

Now consider the last equation, i.e. the equation (4.3).
Difference of the left and the right sides of equation (4.3) we
denote by $K$, i.e.
                 $$
K = \sum\limits_{\gamma  = 0}^N {C[\gamma ]} \frac{{(h\gamma  -
1)^{2m - 2} }}{{2(2m - 2)!}} - \frac{A}{{2(2m - 3)!}} +
\sum\limits_{\alpha  = 1}^{m - 1} {\alpha \lambda _\alpha  }  -
\frac{1}{{2(2m - 1)!}}
$$
and therefore we keep in mind that $K=0$.

Applying binomial formula for first expression of $K$ and taking
into account (4.4)-(4.6), after some simplifications
    $$
K =  - \frac{1}{{2(2m - 1)!}} + \frac{1}{{2(2m - 2)!}} -
\frac{1}{{2(2m - 3)!}}\left( {\frac{1}{2} - B} \right) +
\sum\limits_{j = 1}^{2m - 4} {\frac{{( - 1)^j }}{{2j!(2m - 1 -
j)!}} + }
$$
$$
+ \sum\limits_{j = 1}^{m - 1} {\frac{{( - 1)^j }}{{2(2m - j)!(j -
1)!}} - \sum\limits_{j = m - 1}^{2m - 4} {\frac{{( - 1)^j
B}}{{2j!(2m - 3 - j)!}} - } }
$$
$$
- \sum\limits_{j = 1}^{m - 1} {\frac{{( - 1)^j B}}{{2(2m - 2 -
j)!(j - 1)!}} + \sum\limits_{j = 1}^{m - 1} {\frac{{( - 1)^{j - 1}
}}{{(j - 1)!(2m - 1 - j)!}}\sum\limits_{\gamma  = 0}^N {C[\gamma
](h\gamma )^{2m - 1 - j} }  + } }
$$
$$
+ \sum\limits_{j = 1}^{m - 1} {\frac{{h^{2m - j} }}{{(2m - 1 -
j)!\,(j - 1)!}}\left[ {\sum\limits_{k = 1}^{m - 1} {\sum\limits_{i
= 0}^{2m - 1 - j} {\frac{{d_k q_k  + p_k q_k^{N + i} ( - 1)^{i +
1} }}{{(q_k  - 1)^{i + 1} }}\Delta ^i 0^{2m - 1 - j}  -
\frac{{B_{2m - j} }}{{2m - j}}} } } \right]}.\eqno   (5.31)
$$
Now consider the sum $\sum\limits_{\gamma  = 0}^N {C[\gamma
](h\gamma )^{2m - 1 - j} }$ in (5.31). For this sum using formulas
 (5.1), (2.4)-(2.6) after some simplifications we obtain
$$
\sum\limits_{\gamma  = 0}^N {C[\gamma ](h\gamma )^{2m - 1 - j} }
= \sum\limits_{i = 1}^{2m - j} {\frac{{(2m - 1 - j)!\,B_{2m - j -
i} }}{{i!\,(2m - j - i)!}}h^{2m - j - i}  + }
$$
$$
+ h^{2m - j} \left( \sum\limits_{k = 1}^{m - 1} \sum\limits_{i =
0}^{2m - 1 - j} {\frac{{d_k q_k^i  + p_k q_k^{N + 1} ( - 1)^{i +
1} }}{{(1 - q_k )^{i + 1} }}} \Delta ^i 0^{2m - 1 - j}  - \right.
$$
$$
\left.- \sum\limits_{k = 1}^{m - 1} {\sum\limits_{i = 0}^{2m - 1 -
j} {\frac{{d_k q_k^{N + i}  + p_k q_k^{} ( - 1)^{i + 1} }}{{(1 -
q_k )^{i + 1} }}} \Delta ^i 0^{2m - 1 - j} }   \right) -
$$
$$
- \sum\limits_{p = 0}^{2m - 2 - j} {\frac{{(2m - 1 - j)!\,h^{p +
1} }}{{p!\,(2m - 1 - j - p)!}}\sum\limits_{k = 1}^{m - 1}
{\sum\limits_{i = 0}^p {\frac{{d_k q_k^{N + i}  + p_k q_k^{} ( -
1)^{i + 1} }}{{(1 - q_k )^{i + 1} }}} \Delta ^i 0^p } }. \eqno
(5.32)
$$
Substituting (5.32) into (5.31) we get polynomial of degree $2m -
1$ with respect to $h$. It is easy to see that constant term and
coefficients in front of $h$ and $h^2 $ are zero. Then for $K$ we
obtain
$$
K = \sum\limits_{j = 1}^{m - 1} {\frac{{( - 1)^{j - 1} }}{{(j -
1)!(2m - 1 - j)!}}\left( {\sum\limits_{i = 1}^{2m - 3 - j}
{\frac{{(2m - 1 - j)!B_{2m - j - i} }}{{i!\,(2m - j - i)!}}h^{2m -
j - i}  - } } \right.}
$$
 $$
- \sum\limits_{p = 2}^{2m - 2 - j} {\frac{{(2m - 1 - j)!h^{p + 1}
}}{{p!\,(2m - 1 - j - p)!}}\sum\limits_{k = 1}^{m - 1} {\left.
{\sum\limits_{i = 0}^p {\frac{{d_k q_k^{N + i}  + p_k q_k^{} ( -
1)^{i + 1} }}{{(1 - q_k )^{i + 1} }}} \Delta ^i 0^p } \right)}  +
}
$$
$$
+ \sum\limits_{j = 1}^{m - 1} \frac{{h^{2m - j} }}{{(j - 1)!(2m -
1 - j)!}}\left( {( - 1)^{j - 1} \sum\limits_{k = 1}^{m - 1}
{\sum\limits_{i = 0}^{2m - 1 - j} {\frac{{d_k q_k^i  + p_k q_k^{N
+ 1} ( - 1)^{i + 1} }}{{(1 - q_k )^{i + 1} }}} \Delta ^i 0^{2m - 1
- j}  + } } \right.
$$
$$
+ ( - 1)^j \sum\limits_{k = 1}^{m - 1} \sum\limits_{i = 0}^{2m - 1
- j} {\frac{{d_k q_k^{N + i}  + p_k q_k ( - 1)^{i + 1} }}{{(1 -
q_k )^{i + 1} }}} \Delta ^i 0^{2m - 1 - j}  +
$$
$$
+\sum\limits_{k = 1}^{m - 1} \left. {\sum\limits_{i = 0}^{2m - 1 -
j} {\frac{{d_k q_k  + p_k q_k^{N + i} ( - 1)^{i + 1} }}{{(1 - q_k
)^{i + 1} }}} \Delta ^i 0^{2m - 1 - j} } \right)   -
\sum\limits_{j = 1}^{m - 1} {\frac{{h^{2m - j} B_{2m - j} }}{{(2m
- j)!(j - 1)!}}}.
$$
Hence if take into account lemma 5.2 and (5.4*), then the
expression in the second parenthesis are simplified and $K$ has
the form
    $$
K = \sum\limits_{j = 1}^{m - 1} {\frac{{( - 1)^{j - 1} }}{{(j -
1)!}}} \sum\limits_{i = 1}^{2m - 3 - j} {\frac{{B_{2m - j - i}
h^{2m - j - i} }}{{i!\,(2m - j - i)!}} - } \sum\limits_{j = 1}^{m
- 1} {\frac{{h^{2m - j} B_{2m - j} }}{{(2m - j)!\,(j - 1)!}}}  -
$$
$$
- \sum\limits_{j = 1}^{m - 1} {\frac{{( - 1)^{j - 1} }}{{(j -
1)!}}} \sum\limits_{p = 2}^{2m - 1 - j} {\frac{{h^{p + 1}
}}{{p!\,(2m - 1 - j - p)!}}Z_p }.\eqno   (5.33)
$$
Using lemma 5.3 from (5.33) we get
$$
K = \sum\limits_{j = 3}^m {\frac{{B_j h^j }}{{j!}}} \sum\limits_{i
= 0}^{m - 2} {\frac{{( - 1)^i }}{{i!\,(2m - 1 - j - i)!}} + }
\sum\limits_{j = m + 1}^{2m - 2} {\frac{{B_j h^j }}{{j!}}}
\sum\limits_{i = 0}^{2m - 2 - j} {\frac{{( - 1)^i }}{{i!\,(2m - 1
- j - i)!}} - }
$$
$$
- \sum\limits_{j = 1}^{m - 1} {\frac{{h^{2m - j} B_{2m - j}
}}{{(2m - j)!\,(j - 1)!}}}  - \sum\limits_{j = 3}^{m + 1}
{\frac{{h^j Z_{j - 1} }}{{(j - 1)!}}} \sum\limits_{i = 0}^{m - 2}
{\frac{{( - 1)^i }}{{i!\,(2m - 1 - j - i)!}} - }
$$
$$
- \sum\limits_{j = m + 2}^{2m - 1} {\frac{{h^j Z_{j - 1} }}{{(j -
1)!}}} \sum\limits_{i = 0}^{2m - 1 - j} {\frac{{( - 1)^i
}}{{i!\,(2m - 1 - j - i)!}} = 0}.
$$
Hence using (5.19) after simplifications we have
 $$
K = \frac{{h^m }}{{(m - 1)!}}\left( {\frac{{B_m }}{m} - Z_{m - 1}
} \right)\sum\limits_{i = 0}^{m - 2} {\frac{{( - 1)^i }}{{i!\,\,(m
- 1 - i)!}} + }
$$
$$
+\sum\limits_{j = m + 1}^{2m - 2} {\frac{{B_j h^j
}}{{j!}}\sum\limits_{i = 1}^{2m - 2 - j} {\frac{{( - 1)^i
}}{{i!(2m - 1 - j - i)!}} - } }
$$
$$
- \sum\limits_{j = m + 2}^{2m - 2} {\frac{{h^j Z_{j - 1} }}{{(j -
1)!(2m - 1 - j)!}}\sum\limits_{i = 0}^{2m - 1 - j} {\frac{{(2m - 1
- j)!( - 1)^i }}{{i!\,(2m - 1 - j - i)!}} - \frac{{h^{2m - 1}
Z_{2m - 2} }}{{(2m - 2)!}} = } }
$$
$$
= \frac{{h^m ( - 1)^m }}{{\left[ {(m - 1)!} \right]^2 }}\left(
{\frac{{B_m }}{m} - Z_{m - 1} } \right) + \sum\limits_{j = m +
1}^{2m - 2} {\frac{{B_j h^j }}{{j!\,(2m - 1 - j)!}}\left( {( -
1)^j  - 1} \right) - } \frac{{h^{2m - 1} Z_{2m - 2} }}{{(2m -
2)!}}.
$$
Here  the middle sum is equal to zero because when $j$ is even
 $( - 1)^j - 1 = 0$, and when $j$ is odd Bernoulli numbers
$B_j  = 0$. Therefore finally for  $K$ we get
$$
K =  - \frac{{h^{2m - 1} Z_{2m - 2} }}{{(2m - 2)!}} + \frac{{h^m (
- 1)^m }}{{\left[ {(m - 1)!} \right]^2 }}\left( {\frac{{B_m }}{m}
- Z_{m - 1} } \right) = 0.
$$
It means hence taking into account (5.4*) we get following
equation
$$
\sum\limits_{k = 1}^{m - 1} {\sum\limits_{i = 0}^{m - 1}
{\frac{{d_k q_k^{N + i}  + p_k q_k^{} ( - 1)^{i + 1} }}{{(1 - q_k
)^{i + 1} }}} \Delta ^i 0^{m - 1}  = \frac{{B_m }}{m},} \eqno
(5.34)
$$
$$
\sum\limits_{k = 1}^{m - 1} {\sum\limits_{i = 0}^{2m - 2}
{\frac{{d_k q_k^{N + i}  + p_k q_k^{} ( - 1)^{i + 1} }}{{(1 - q_k
)^{i + 1} }}} \Delta ^i 0^{2m - 2}  = 0}. \eqno
(5.35)
$$
Using (5.30) when $\alpha  = 1$ and taking into account lemma 5.2
it is easy to show , that the equation  (5.34) coincide by
equation of the system (5.19) when $\alpha  = 1$. Thus, we have
obtained the last equation (5.35) for $d_k $ and $p_k $, which is
the same as (5.13). Theorem is proved.

\subsection{Calculation of norm of the error functional}

For square of norm of the error functional (1.2) of optimal
quadrature formulas of the form (1.1) take placed following
theorem.

\textbf{Theorem 5.2.} \textit{For square of norm of the error
functional  (1.2) of optimal quadrature formula of the form (1.1)
following is valid
    $$
\left\| \ell|L_2^{(m)*}(0,1)\right\|^2  = ( - 1)^{m + 1} \left[
{\frac{{B_{2m} h^{2m} }}{{(2m)!}}} \right. + \frac{{h^{2m + 1}
}}{{(2m)!}}\left. {\sum\limits_{k = 1}^{m - 1} {\sum\limits_{i =
0}^{2m} {\frac{{(1 - q_k^N )(d_k q_k^i  + ( - 1)^i p_k q_k )}}{{(1
- q_k )^{i + 1} }}} } \Delta ^i 0^{2m} } \right],
$$
where $d_{k,} \,p_k $ are determined from system (5.10)-(5.13),
$B_{2m} $ are Bernoulli numbers, $q_k $ are roots of Euler
polynomial of degree $2m - 2$, $\,\left| {q_k } \right| < 1.$ }

\textbf{Proof.} Computing defined integrals in the expression
$||\ell ||^2 $ we get
$$
||\ell ||^2  = ( - 1)^{m + 1} \left[ {\frac{{A \cdot B}}{{(2m -
3)!}} - \frac{{A - B}}{{(2m - 1)!}} + \sum\limits_{\beta  = 0}^N
{C[\beta ]} } \right.\frac{{A\,(h\beta )^{2m - 2}  - B\,(h\beta  -
1)^{2m - 2} }}{{2(2m - 2)!}} +
$$
$$
+ \sum\limits_{\beta  = 0}^N {C[\beta ]} \,F(h\beta ) +
\sum\limits_{\beta  = 0}^N {C[\beta ]} \left\{ {F(h\beta ) -
\sum\limits_{\gamma  = 0}^N {C[\gamma ]\frac{{|h\beta  - h\gamma
|^{2m - 1} }}{{2(2m - 1)!}}}  + } \right.
$$
$$
+ \left. {\frac{{A\,(h\beta )^{2m - 2}  - B\,(h\beta  - 1)^{2m -
2} }}{{2(2m - 2)!}}} \right\} - \left. {\frac{1}{{(2m + 1)!}}}
\right],
$$
where $F(h\beta )$ is determine by formula (5.16).  As is obvious
from here according to (4.1) the expression into curly brackets is
equal to the polynomial $P_{m - 1} (h\beta ) = \sum\limits_{\alpha
= 0}^{m - 1} {\lambda _\alpha (h\beta )^\alpha } $. Then $||\ell
||^2 $ have the form
$$
||\ell ||^2  = ( - 1)^{m + 1} \left[ {\frac{{A \cdot B}}{{(2m -
3)!}} - \frac{{A - B}}{{(2m - 1)!}} + \sum\limits_{\beta  = 0}^N
{C[\beta ]} } \right.\frac{{A\,(h\beta )^{2m - 2}  - B\,(h\beta  -
1)^{2m - 2} }}{{2(2m - 2)!}} +
$$
$$
+ \sum\limits_{\beta  = 0}^N {C[\beta ]} \,F(h\beta ) +
\sum\limits_{\beta  = 0}^N {C[\beta ]} \,\,P_{m - 1} (h\beta ) -
\left. {\frac{1}{{(2m + 1)!}}} \right].
$$
Hence using (4.2) and (4.3) we get
$$
||\ell ||^2  = ( - 1)^{m + 1} \left[ {\frac{{A \cdot B}}{{(2m -
3)!}} - \frac{{A - B}}{{(2m - 1)!}} + } \right.A \cdot \left(
{\frac{1}{{2(2m - 1)!}} + \lambda _1  - \frac{B}{{2(2m - 3)!}}}
\right) -
$$
$$
- B \cdot \left( {\frac{1}{{2(2m - 1)!}} - \sum\limits_{\alpha  =
1}^{m - 1} {\alpha \lambda _\alpha   + \frac{A}{{2(2m - 3)!}}} }
\right) +
$$
$$
+ \sum\limits_{\beta  = 0}^N {C[\beta ]} \,F(h\beta ) +
\sum\limits_{\beta  = 0}^N {C[\beta ]} \,\,P_{m - 1} (h\beta ) -
\left. {\frac{1}{{(2m + 1)!}}} \right].        \eqno   (5.35)
$$
From (5.35) after simplifications using (5.17), (5.12) we have
$$
||\ell ||^2  = ( - 1)^{m + 1} \left[ {\frac{{B - A}}{{2(2m - 1)!}}
+ } \right.\frac{A}{{(2m - 2)!}} \cdot \left( { - \frac{1}{{2(2m -
1)}} + \frac{1}{2}\sum\limits_{\gamma  = 0}^N {C[\gamma ]}
(h\gamma )^{2m - 2} } \right. +
$$
$$
+ \left. {\frac{{(2m - 2)B}}{2}} \right) + B \cdot \sum\limits_{j
= 1}^{m - 1} {\frac{1}{{(2m - 1 - j)!(j - 1)!}}\left( {\frac{{( -
1)^{2m - j} }}{{2(2m - j)}} - \frac{{B_{2m - j} h^{2m - j} }}{{2m
- j}} - } \right.}
$$
$$
- h^{2m - j} \sum\limits_{k = 1}^{m - 1} {\sum\limits_{i = 0}^{2m
- 1 - j} {\frac{{ - d_k q_k  + p_k q_k^{N + i} ( - 1)^i }}{{(q_k
- 1)^{i + 1} }}\Delta ^i 0^{2m - 1 - j}  +
\frac{1}{2}\sum\limits_{\gamma  = 0}^N {C[\gamma ]( - h\gamma
)^{2m - 1 - j}  - } } }
$$
$$
\left. { - (2m - 1 - j)B\frac{{( - 1)^j }}{2}} \right) +
\sum\limits_{\beta  = 0}^N {C[\beta ]\,F(h\beta ) + \left.
{\sum\limits_{\beta  = 0}^N {C[\beta ]\,P_{m - 1} (h\beta ) -
\frac{1}{{(2m + 1)!}}} } \right]}.
$$

Hence, taking into account (5.16), (5.17), (5.18) and using
(4.4)-(4.6), after some calculations we obtain
$$
||\ell ||^2  = ( - 1)^{m + 1} \left[ {\frac{B}{{(2m - 1)!}} -
\frac{1}{{(2m + 1)!}} - } \right.\sum\limits_{j = 0}^{m - 3}
{\frac{{B\,\,( - 1)^j }}{{(2m - 2 - j)!\,(j + 1)!}} +
\sum\limits_{j = 0}^{m - 1} {\frac{{( - 1)^j }}{{(j + 1)!\,(2m -
j)!}} - } }
$$
$$
- \sum\limits_{j = 2}^{m - 1} {\frac{{h^{2m - j} }}{{(2m - 1 -
j)!(j + 1)!}}\left( {\frac{{B_{2m - j} }}{{2m - j}} +
\sum\limits_{k = 1}^{m - 1} {\sum\limits_{i = 0}^{2m - 1 - j}
{\frac{{ - d_k q_k  + p_k q_k^{N + i} ( - 1)^i }}{{(q_k  - 1)^{i +
1} }}\Delta ^i 0^{2m - 1 - j} } } } \right)}  +
$$
$$
+ \sum\limits_{j = 0}^{m - 1} {\frac{{( - 1)^{2m - 1 - j} }}{{(2m
- 1 - j)!(j + 1)!}}\sum\limits_{\gamma  = 0}^N {C[\gamma ](h\gamma
)^{2m - 1 - j}  + } \left. {\sum\limits_{\beta  = 0}^N {C[\gamma
]\frac{{(h\beta )^{2m} }}{{(2m)!}}} } \right]}.\eqno   (5.36)
$$
When $\alpha  > m - 1$ using lemma 5.1 and formulas (2.4)-(2.6) we
get
$$
\sum\limits_{\gamma  = 0}^N {C[\gamma ](h\gamma )^\alpha   =
\frac{1}{{\alpha  + 1}} + \sum\limits_{j = 1}^{\alpha  - 1}
{\frac{{\alpha !B_{\alpha  + 1 - j} }}{{j!(\alpha  + 1 -
j)!}}h^{\alpha  + 1 - j}  + } }
$$
$$
+ h^{\alpha  + 1} \sum\limits_{k = 1}^{m - 1} {\sum\limits_{i =
0}^\alpha  {\frac{{d_k q_k^i  + p_k q_k^{N + 1} ( - 1)^{i + 1}
}}{{(1 - q_k )^{i + 1} }}\Delta ^i 0^\alpha   - } }
$$
$$
- \sum\limits_{j = 1}^\alpha  {\frac{{\alpha !h^{j + 1}
}}{{j!(\alpha  - j)!}}\sum\limits_{k = 1}^{m - 1} {\sum\limits_{i
= 0}^j {\frac{{d_k q_k^{N + i}  + p_k q_k^{} ( - 1)^{i + 1} }}{{(1
- q_k )^{i + 1} }}\Delta ^i 0^j } } }.\eqno   (5.37) $$ Using
(5.37) and taking into account (2.1), (2.3), after simplications,
from (5.36) we have
$$
||\ell ||^2  = ( - 1)^{m + 1} \left[ {\sum\limits_{j = m}^{2m - 3}
{\frac{{B_{2m - j} h^{2m - j} }}{{(2m - j)!(j + 1)!}} +
\frac{{B_{2m} h^{2m} }}{{(2m)!}} + } } \right.
$$
$$
+ \frac{{h^{2m + 1} }}{{(2m)!}}\sum\limits_{k = 1}^{m - 1}
{\sum\limits_{i = 0}^{2m} {\frac{{d_k q_k^i  + p_k q_k^{N + 1} ( -
1)^{i + 1} }}{{(1 - q_k )^{i + 1} }}\Delta ^i 0^{2m}  + } }
$$ $$
+\sum\limits_{\alpha  = m}^{2m - 1} {\frac{{( - 1)^\alpha  }}{{(2m
- \alpha )!}}\sum\limits_{j = 1}^{\alpha  - 2} {\frac{{B_{\alpha +
1 - j} h^{\alpha  + 1 - j} }}{{j!(\alpha  + 1 - j)!}} - } }
$$
$$
- \sum\limits_{\alpha  = m}^{2m} {\frac{{( - 1)^\alpha  }}{{(2m -
\alpha )!}}\sum\limits_{j = 2}^\alpha  {\frac{{h^{j + 1}
}}{{j!(\alpha  - j)!}}\left. {\sum\limits_{k = 1}^{m - 1}
{\sum\limits_{i = 0}^j {\frac{{d_k q_k^{N + i}  + p_k q_k^{} ( -
1)^{i + 1} }}{{(1 - q_k )^{i + 1} }}\Delta ^i 0^j } } } \right].}
}
$$
Hence last two sums regrouping in powers of $h$, using designation
 (5.4*) and keepin in mind $Z_{\alpha  - 1} =
\frac{{B_\alpha }}{\alpha },\,\,\,\alpha  = \overline {3,m} $, we
have
$$
||\ell ||^2  = ( - 1)^{m + 1} \left[ {\sum\limits_{j = m}^{2m - 3}
{\frac{{B_{2m - j} h^{2m - j} }}{{(2m - j)!(j + 1)!}} +
\frac{{B_{2m} h^{2m} }}{{(2m)!}} + } } \right.
$$
$$
+ \frac{{h^{2m + 1} }}{{(2m)!}}\sum\limits_{k = 1}^{m - 1}
{\sum\limits_{i = 0}^{2m} {\frac{{d_k q_k^i  + p_k q_k^{N + 1} ( -
1)^{i + 1} }}{{(1 - q_k )^{i + 1} }}\Delta ^i 0^{2m}  + } }
$$
$$
+ \sum\limits_{\alpha  = 3}^m {\frac{{B_\alpha  h^\alpha
}}{{\alpha !}}\sum\limits_{j = m}^{2m - 1} {\frac{{( - 1)^j
}}{{(2m - j)!\,\,(j - \alpha  + 1)!}} + } } \sum\limits_{\alpha  =
m + 1}^{2m - 1} {\frac{{B_\alpha  h^\alpha  }}{{\alpha
!}}\sum\limits_{j = \alpha }^{2m - 1} {\frac{{( - 1)^j }}{{(2m -
j)!\,\,(j - \alpha  + 1)!}} - } }
$$
$$
- \sum\limits_{\alpha  = 3}^m {\frac{{B_\alpha  h^\alpha
}}{{\alpha !}}\sum\limits_{j = m}^{2m} {\frac{{( - 1)^j }}{{(2m -
j)!\,\,(j - \alpha  + 1)!}} - } } \sum\limits_{\alpha  = m +
1}^{2m + 1} {\left. {\frac{{Z_{\alpha  - 1} h^\alpha  }}{{(\alpha
- 1)!}}\sum\limits_{j = \alpha  - 1}^{2m} {\frac{{( - 1)^j }}{{(2m
- j)!\,\,(j - \alpha  + 1)!}}} } \right]}.\eqno   (5.38)
$$
Since in (5.38)
$$
\sum\limits_{j = \alpha  - 1}^{2m} {\frac{{( - 1)^j }}{{(2m -
j)!\,(j - \alpha  + 1)!}} = \frac{{( - 1)^{\alpha  - 1} }}{{(2m -
\alpha  - 1)!}}(1 - 1)^{2m - \alpha  + 1}  = 0,\,} \mbox{ при  } m
+ 1 \le \alpha  \leq 2m
$$
and
$$
\sum\limits_{j = \alpha }^{2m - 1} {\frac{{( - 1)^j }}{{(2m -
j)!\,\,(j - \alpha  + 1)!}} =  - \frac{1}{{(2m - \alpha  + 1)!}}((
- 1)^{\alpha  - 1}  + 1)},
$$
then, using these equalities,  from (5.38) we obtain
$$
\left\| \ell  \right\|^2  = ( - 1)^{m + 1} \left[ {\frac{{B_{2m}
h^{2m} }}{{(2m)!}}} \right. + \frac{{h^{2m + 1}
}}{{(2m)!}}\sum\limits_{k = 1}^{m - 1} {\sum\limits_{i = 0}^{2m}
{\frac{{(1 - q_k^N )(d_k q_k^i  + ( - 1)^i p_k q_k )}}{{(1 - q_k
)^{i + 1} }}} } \Delta ^i 0^{2m}  -
$$
$$
- \left. {\sum\limits_{j = m + 1}^{2m - 1} {\frac{{B_\alpha \,\,((
- 1)^{\alpha  - 1}  + 1)}}{{\alpha !\,\,(2m - \alpha  + 1)!}}} }
\right]. $$ Hence taking into account that when  $\alpha $ is even
$ ( - 1)^{\alpha - 1} + 1 = 0$ and when  $\alpha $ is odd
$B_\alpha = 0$, (since $\alpha \ne 1$), we get the statement of
theorem 5.2.
Theorem 5.2 is proved.

\textbf{Corollary 5.1.} \textit{In the space $L_2^{(2)}(0,1)$
among quadrature formulas of the form (1.1) with the error
functional (1.2) there exists unique optimal formula  which
coefficients are determined by following formulas}
    $$
C[\beta] = \left\{
\begin{array}{ll}
\frac{h}{2},& \beta=0,N,\\
h,&\beta=\overline{1,N-1},\\
\end{array} \right.
$$
$$
A =\frac{h^2}{12},\ \ \  B =-\frac{h^2}{12}.
$$
\textit{Furthermore for square of norm of the error functional
following is valid}
$$
\left\| \ell|L_2^{(2)*}(0,1) \right\|^2  = \frac{h^4}{720}.
$$

\textbf{Corollary 5.2.} \textit{In the space $L_2^{(3)}(0,1)$
among quadrature formulas of the form (1.1) with the error
functional (1.2) there exists unique optimal formula  which
coefficients are determined by following formulas }
    $$
C[\beta] = \left\{
\begin{array}{ll}
\frac{h}{2},& \beta=0,N,\\
h,&\beta=\overline{1,N-1},\\
\end{array} \right.
$$
$$
A =\frac{h^2}{12},\ \ \  B =-\frac{h^2}{12}.
$$
\textit{Furthermore for square of norm of the error functional
following is valid}
$$
\left\| \ell|L_2^{(3)*}(0,1) \right\|^2  = \frac{h^6}{30240}.
$$

Proofs of Corollaries 5.1 and 5.2 we get immediately from theorems
5.1 and 5.2 when $m=2$ and $m=3$ respectively.

\section{Acknowledgements}

The second author gratefully acknowledges the Abdus Salam School
of Mathematical Sciences (ASSMS), GC University, Lahore, Pakistan
for
providing the Post Doctoral Research Fellowship.\\[0.2cm]

\textbf{Shadimetov Kholmat Mahkambaevich}\\
Institute of Mathematics and Information Technologies\\
Uzbek Academy of Sciences\\
29, Durmon yuli str \\
 Tashkent, 100125, Uzbekistan

\textbf{Hayotov Abdullo Rahmonovich}\\
Institute of Mathematics and Information Technologies\\
Uzbek Academy of Sciences\\
29, Durmon yuli str \\
 Tashkent, 100125, Uzbekistan\\
e-mail: hayotov@mail.ru, abdullo\_hayotov@mail.ru

\textbf{Nuraliev Farhod Abduganievich}\\
Institute of Mathematics and Information Technologies\\
Uzbek Academy of Sciences\\
29, Durmon yuli str \\
Tashkent, 100125, Uzbekistan

\end{document}